\documentclass[11pt,reqno]{amsart}
\usepackage[latin1]{inputenc}
\usepackage{amsmath}
\usepackage{amsfonts}
\usepackage{amssymb,esint}

\theoremstyle{plain}
\newtheorem*{theorem}{Theorem}
\newtheorem*{lemma}{Lemma}

\newcommand{\R}{\mathbb R}
\newcommand{\N}{\mathbb N}

\usepackage{color}

\newcommand{\res}
{\mathop{\hbox{\vrule height 7pt width .5pt depth 0pt \vrule
height .5pt width 6pt depth 0pt}}\nolimits}

\long\def\MSC#1\EndMSC{\def\arg{#1}\ifx\arg\empty\relax\else
     {\par\narrower\noindent
     {\small\it 2010 Mathematics Subject Classification.} \small #1\par}\fi}

\long\def\KEY#1\EndKEY{\def\arg{#1}\ifx\arg\empty\relax\else
     {\par\narrower\noindent
     {\small\it Keywords and Phrases.} \small #1\par}\fi}

\author[Massaccesi]{Annalisa Massaccesi}
\author[Vittone]{Davide Vittone}
\address[Massaccesi]{ Institut f\"ur Mathematik, Winterthurerstrasse 190, CH-8057 Z\"urich, Switzerland.}
\email{annalisa.massaccesi@math.uzh.ch }
\address[Vittone]{Dipartimento di Matematica, via Trieste 63, 35121 Padova, Italy.}
\email{vittone@math.unipd.it}
\thanks{A.M. is supported by Universit\"at Z\"urich and ERC grant RAM (Regularity for Area Minimizing currents), ERC 306247. D.V. is supported by University of Padova and GNAMPA of INDAM (Italy). He also wishes to thank the Institut f\"ur Mathematik, Z\"urich, for its hospitality during the preparation of this paper.}

\begin{document}


\title{An elementary proof of the rank-one theorem for BV functions}

\begin{abstract}
We provide a simple proof of a result, due to G. Alberti, concerning a rank-one property for the singular part of the derivative of vector-valued functions of bounded variation.
\end{abstract}

\maketitle

In this paper we provide a short, elementary proof of the following result by G. Alberti \cite{Alb-RankOne} concerning a rank-one property for the derivative of a function with bounded variation.

\begin{theorem}
Let $\Omega\subset\R^n$ be an open set, $u:\Omega\to\R^m$  a function of bounded variation and let $D_su$ be the singular part of $Du$ with respect to the Lebesgue measure $\mathcal L^n$. Then $D_su$ is a rank-one measure, i.e., the (matrix-valued) function $\frac{D_su}{|D_su|}(x)$ has rank one for $|D_su|$-a.e. $x\in\Omega$.
\end{theorem}

We recall that a function $u\in L^1(\Omega,\R^m)$ has {\em bounded variation} in $\Omega$ ($u\in BV(\Omega,\R^m)$) if the derivatives $Du$ of $u$ in the sense of distributions are represented by a (matrix-valued) measure with finite total variation. The measure $Du$ can then be decomposed as the sum $Du=D_au+D_su$ of a measure $D_au$, that is absolutely continuous with respect to $\mathcal L^n$, and a measure $D_su$ that is singular with respect to $\mathcal L^n$. The Radon-Nikodym derivative $\frac{D_su}{|D_su|}$ of $D_su$ with respect to its total variation $|D_su|$ is a $|D_su|$-measurable map from $\Omega$ to $\R^{m\times n}$. The Theorem  states that this map takes values in the space of rank-one matrices. See \cite{AFP} for more details on $BV$ functions.

The Theorem above was conjectured by L. Ambrosio and E. De Giorgi in \cite{ADGNuovoTipoFunz}. It was first proved by G. Alberti in \cite{Alb-RankOne} by introducing new tools and using  sophisticated techniques in Geometric Measure Theory. A new proof has been announced by G. De Philippis and F. Rindler as one of the consequences of the forthcoming, profound PDE result  \cite{DP}.

On the contrary, our proof of the Theorem above is elementary: it stems from well-known geometric properties relating the derivative of a BV function and the \emph{perimeter} of its subgraph. The main new tool is the following lemma, where we denote by $\pi:\R^{n+1}\to\R^n$ the canonical projection $\pi(x_1,\ldots,x_{n+1}):=(x_1,\ldots, x_n)$. 


\begin{lemma}
Let $\Sigma_1,\Sigma_2$ be $C^1$  hypersurfaces in $\R^{n+1}$ with unit normals $\nu_{\Sigma_1},\nu_{\Sigma_2}$. 
Then, the set
\[
T:=\left\{p\in \Sigma_1: 
\exists\,q\in\Sigma_2\cap\pi^{-1}(\pi(p))\text{ with }\left(\nu_{\Sigma_1}(p)\right)_{n+1}=\left(\nu_{\Sigma_2}(q)\right)_{n+1}=0
\text{ and }
\nu_{\Sigma_1}(p)\neq\pm\nu_{\Sigma_2}(q)
\right\}
\]
is $\mathcal H^n$-negligible.
\end{lemma}

We postpone the proof of the Lemma  in order to directly address the proof of the main result.

\begin{proof}[Proof of the Theorem]
Let $u=(u_1,\dots,u_m)\in BV(\Omega,\R^m)$. It is not restrictive to assume that $\Omega$ is bounded. For any $i=1,\dots,m$ we write $D_su_i=\sigma_i|D_su_i|$ for a $|D_su_i|$-measurable map $\sigma_i:\Omega\to {\mathbb S}^{n-1}$. We also let $E_i:=\{(x,t)\in\Omega\times\R:\,t<u_i(x)\}$ be the subgraph of $u_i$; it is well known that $E_i$ has {finite perimeter} in $\Omega\times\R$. Denoting by $\partial^* E_i$ the \emph{reduced boundary} of $E_i$ and by $\nu_i$ the measure theoretic inner normal to $E_i$, we have (see e.g. \cite[Section 4.1.5]{GMS})
\[
|D_su_i|=\pi_\#({\mathcal H}^{n}\res S_i)\quad\text{for }S_i:=\left\{p\in\partial^\ast E_i:\left(\nu_i(p)\right)_{n+1}=0\right\},
\]
where $\pi_\#$ denotes push-forward of measures. The set $S_i$ is $n$-rectifiable and we can assume that it is contained in the union $\cup_{h\in\N} \Sigma^i_h$ of $C^1$ hypersurfaces $\Sigma^i_h$ in $\R^{n+1}$. 

By  \cite[Section 4.1.5]{GMS}, the Lemma above and well-known properties of rectifiable sets, the following properties hold for $\mathcal H^n$-a.e. $p\in S_i$:
\begin{align}
& \nu_{\partial^\ast E_i}(p) = (\sigma_i(\pi(p)),0)  \label{2.1}\\
& \text{if }p\in\Sigma^i_h, \text{ then }  \nu_{i}(p) = \pm\nu_{\Sigma^i_h}(p) \label{2.2}\\
& \text{if }p\in \Sigma^i_h\text{ and }  q\in S_j\cap\Sigma^j_k\cap\pi^{-1}(\pi(p)),\text{ then }  \nu_{\Sigma^i_h}(p)=\pm \nu_{\Sigma^j_k}(q). \label{2.3}
\end{align}
Up to modifying $S_i$ on a ${\mathcal H}^n$-negligible set and $\sigma_i$ on a $|D_su_i|$-negligible set, we can assume that  \eqref{2.1},\eqref{2.2} and \eqref{2.3} hold everywhere on $S_i$ and that $\sigma_i=0$ on $\Omega\setminus\pi(S_i)$.

Since $D_su=(\sigma_1|D_su_1|,\dots,\sigma_m|D_su_m|)$ and $|D_su|$ is 
concentrated on $\pi(S_1)\cup\dots\cup \pi(S_m)$, it is enough to prove that the matrix-valued function $(\sigma_1,\dots,\sigma_m)$ has rank 1 on $\pi(S_1)\cup\dots\cup \pi(S_m)$. This will follow if we prove that the implication
\[
i, j\in\{1,\ldots,m\},\ i\neq j,\ x\in\pi(S_i)\ \Longrightarrow\ \sigma_j(x)\in\{0,\sigma_i(x),-\sigma_i(x)\}
\]
holds. If $i,j,x$ are as above and $x\notin \pi(S_j)$, then $\sigma_j(x)=0$. Otherwise, $x\in \pi(S_i)\cap \pi(S_j)$, i.e., there exist $p\in S_i$ and $h\in\N$ such that $\pi(p)=x$ and $\sigma_i(x)=\pm\nu_{\Sigma^i_h}(p)$ and there exist $q\in S_j$ and $k\in\N$ such that $\pi(q)=x$ and $\sigma_j(x)=\pm\nu_{\Sigma^j_k}(p)$. By \eqref{2.3} we obtain $\sigma_j(x)=\pm \sigma_i(x)$, as wished.
\end{proof}

\begin{proof}[Proof of the Lemma]
Consider the sets $\Sigma:=\Sigma_1\times \Sigma_2\subset \R^{2n+2}$ and 
\[
\Delta:=\{\xi=(x,t,y,s)\in\R^n\times\R\times\R^n\times\R=\R^{2n+2}:\,x=y\}.
\]
Then $\Sigma$ is a $2n$-dimensional manifold of class $C^1$ and $\Delta$ is a smooth $(n+2)$-dimensional manifold in $\R^{2n+2}$. Let us consider the set
\[
R:=\{\xi=(x,t,x,s)\in \Delta\cap\Sigma: (\nu_{\Sigma_1}(x,t))_{n+1}=(\nu_{\Sigma_2}(x,s))_{n+1}=0\text{ and } \nu_{\Sigma_1}(x,t)\neq \pm\nu_{\Sigma_2}(x,s)\}.
\]
By construction, the intersection between $\Delta$ and $\Sigma$ is transversal at every point of $R$, thus $R$ is contained in a $n$-dimensional submanifold $\tilde R\subset \Delta\cap\Sigma$ of class $C^1$.  
Let $\phi$ be the projection  $\phi(x,t,y,s):=(x,t)$; notice that $\phi(\tilde R)\subset\Sigma_1$ and $T=\phi(R)$. Moreover, for every  $\xi\in R$ the differential $d\phi_\xi:T_\xi \tilde R\to T_{\phi(\xi)}\Sigma_1$ is not surjective, because the vector $(0,\ldots,0,1)$ is in the kernel of $d\phi_\xi$. By the area formula we deduce that $\mathcal H^n(T)=\mathcal H^n(\phi(R))=0$, as desired.
\end{proof}

\noindent {\em Acknowledgements.} The authors are grateful to G. Alberti for several suggestions and fruitful discussions.

\bibliographystyle{acm}

\end{document}